\documentclass[journal]{IEEEtran}
\usepackage{amsmath}
\usepackage{amsfonts}
\usepackage{graphicx}
\usepackage{epstopdf}
\usepackage{color}
\usepackage{subeqnarray}
\usepackage{algorithm}
\usepackage{algpseudocode}
\usepackage{indentfirst} 
\usepackage{tikz}
\usepackage{adjustbox}
\usepackage[caption=false]{subfig}

\usepackage[utf8]{inputenc}
\usepackage[english]{babel}
\usepackage[T1]{fontenc}

\usepackage{multirow}

\newcommand{\norm}[1]{\lVert#1\rVert}

\algnewcommand{\Initialize}[1]{
  \State \textbf{Initialize:}
  \Statex \hspace*{\algorithmicindent}\parbox[t]{.95\linewidth}{\raggedright #1}
}

\usetikzlibrary{shapes.geometric, arrows}
\tikzstyle{block} = [rectangle, rounded corners, minimum width=3cm, minimum height=1cm, text centered, draw=blue, fill=red!40]
\tikzstyle{joint} = [circle, minimum width=0.1cm, minimum height=0.1cm, text centered, draw=blue, text centered, fill=red!40]
\tikzstyle{arrow} = [thick,->,>=stealth]

\begin{document}
\title{An Efficient Move Blocking Strategy for Multiple Shooting based Nonlinear Model Predictive Control}


\author{Yutao~Chen$^{1}$ \thanks{$^{1}$Yutao~Chen is with the Department of Electrical Engineering, Eindhoven University of Technology, The Netherlands. Email: y.chen2@tue.nl},Nicolò Scarabottolo$^{2}$, Mattia Bruschetta$^{2}$, and Alessandro~Beghi%
\thanks{$^{2}$Nicolò Scarabottolo, Mattia Bruschetta, and Alessandro~Beghi are with the Department of Information Engineering, University of Padova, Italy.} }

\maketitle

\begin{abstract}
	Move blocking (MB) is a widely used strategy to reduce the degrees of freedom of the Optimal Control Problem (OCP) arising in receding horizon control. The size of the OCP is reduced by forcing the input variables to be constant over multiple discretization steps. In this paper, we focus on developing computationally efficient MB schemes for multiple shooting based nonlinear model predictive control (NMPC). The degrees of freedom of the OCP is reduced by introducing MB in the shooting step, resulting in a smaller but sparse OCP. Therefore, the discretization accuracy and level of sparsity is maintained. A condensing algorithm that exploits the sparsity structure of the OCP is proposed, that allows to reduce the computation complexity of condensing from quadratic to linear in the number of discretization nodes. As a result, active-set methods with warm-start strategy can be efficiently employed, thus allowing the use of a longer prediction horizon. A detailed comparison between the proposed scheme and the nonuniform grid NMPC is given. Effectiveness of the algorithm in reducing computational burden while maintaining optimization accuracy and constraints fulfillment is shown by means of simulations with two different problems.
\end{abstract}

\section{Introduction}
With the fast increase of real-time Nonlinear Model Predictive Control (NMPC) applications, research efforts have been focused on finding efficient solutions of on-line optimization problems. The computational burden of NMPC strongly depends on the dimension of the Optimal Control Problem (OCP) that must be recursively solved on-line. The dimension of the OCP increases proportionally with the length of the prediction horizon, the number of discretization nodes, and the dimension of state and control spaces. 

Since a sufficiently long prediction horizon is essential to guarantee stability of NMPC, several methods have been proposed to reduce the number of discretization nodes and consequently the dimension of the OCP. In particular, in the so-called nonuniform grid schemes, non-equidistant discretization nodes along the prediction horizon are adopted, typically more dense at the beginning of the prediction horizon and more sparse at the end of horizon \cite{quirynen2015multiple}. However, the reduction of the number of nodes comes at the cost of a loss of discretization accuracy. To deal with this issue, the number of discretization intervals has been determined a priori off-line while guaranteeing an upper limit of the discretization error in \cite{lazutkin2018approach}. In addition, adaptive time-mesh refinement techniques have been proposed to improve computation efficiency while maintaining a certain degree of discretization accuracy \cite{paiva2015adaptive, potena2018non, lee2018mesh}. The distribution of discretization nodes is adjusted on-line according to specially designed rules and the total number of nodes is bounded. 
Although discretization accuracy is improved, such methods result in time-varying computational time at each sampling instant due to the time-varying number of discretization nodes.

Another popular class of methods aiming at reducing the number of discretization nodes is input move blocking (MB), that reduces the degrees of freedom (DoFs) of the OCP by constraining the input to be constant over several discretization time steps within the prediction horizon. MB has been widely studied and applied in linear MPC problems. A survey of common MB strategies has been given in \cite{cagienard2007move}. Thorough theoretical analyses have been given in \cite{gondhalekar2009controlled, gondhalekar2010least, shekhar2015optimal} with emphasis on properly choosing the block structure to ensure recursive feasibility and to maintain or maximize the region of attraction (ROA).

Nevertheless, MB for NMPC has not yet been thoroughly studied. Theoretical guidelines for the choice of the block structure for MB NMPC are still under development. Relevant studies can be found for non-uniform grid schemes to guarantee recursive feasibility and closed-loop stability \cite{yu2016stable, paiva2018sufficient}. However, computationally efficient algorithms for MB NMPC are not addressed. When MB is applied, existing algorithms benefit only from reduced number of decision variables but do not exploit the structure of the problem.

In this paper, we focus on computationally efficient numerical algorithms for MB NMPC. 
Our main contribution is a sparsity-preserving input MB scheme for multiple shooting based NMPC for real-time applications. The goal is to reduce the on-line computational burden when using a long prediction horizon, while preserving to a large extent numerical and control performance. In the proposed scheme, MB is introduced in the shooting step when discretizing the OCP to a Nonlinear Programming (NLP) problem using multiple shooting \cite{bock1984multiple}. The NLP problem with reduced DoFs is then solved by the Sequential Quadratic Programming (SQP) method. The resulting Quadratic Programming (QP) problem has reduced DoFs but maintains the level of sparsity, thus maintaining the degree of optimality of the solution and rate of convergence of the optimization algorithms \cite{quirynen2015multiple}. To further reduce problem dimension, a tailored condensing algorithm is developed to exploits the sparsity structure: the computational complexity of condensing is reduced from quadratic to linear in the number of discretization nodes. As a result, active-set methods, that typically require a condensing step, can be efficiently employed given a long prediction horizon, taking advantage of on-line warm-start strategies. In this paper, the proposed strategy is applied to the Real-Time Iteration (RTI) scheme, an effective and well known sub-optimal NMPC algorithm \cite{diehl2002real}, where only one QP problem is solved and a sub-optimal solution is obtained at each sampling instant. Note that methods for designing the block structure to satisfy specific control requirements, e.g. maximizing ROA, are not addressed.

To demonstrate the proposed MB strategy, two different problems are considered, namely, the control of an inverted pendulum and the design of a motion cueing strategy for a dynamic driving simulator. Control, numerical and computational performance of the proposed method are compared against that of 1) the standard unblocked NMPC and 2) NMPC with non-uniform discretization grid. Simulation results show that the proposed algorithm can significantly reduce the on-line computation time while maintaining closed-loop performance, without the need of either changing tuning parameters of NMPC or introducing additional time-mesh refinement steps. In addition, the loss of solution optimality and constraints fulfillment due to reduced DoFs is considerably mitigated thanks to sparsity preservation.

The paper is organized as follows. In Section \ref{sec2} the problem formulation based on multiple shooting is presented. MB is introduced in Section \ref{sec3} with a detailed description of the proposed algorithm and analyses on convergence and stability. In Section \ref{sec4} the two application examples are considered to analyze the performance of the proposed approach. Finally, conclusions are drawn in Section \ref{sec5}.

\section{Problem formulation}\label{sec2}
In NMPC, the OCP that has to be solved at each sampling instant has the following general form:
\begin{equation}
\begin{aligned}
\min_{x(\cdot),u(\cdot)} \quad& J= \int_{t_0}^{t_f} h(t,x(t),u(t);p)\,dt+h_f(x(t_f))\\
s.t. \quad & 0 = \hat{x}_0 - x(t_0)\\
& 0 = f(t,\dot{x}(t), x(t),u(t);p), \quad \forall t\in \left [t_0,t_f\right ],\\
& 0 \geq r(x(t),u(t);p), \quad \forall t\in \left [t_0,t_f\right ],\\
& 0 \geq l(x(t_f);p),
\end{aligned}
\label{OCP}
\end{equation}
where $x\in \mathbb{R}^{n_x}$, $u \in \mathbb{R}^{n_u}$, $p\in \mathbb{R}^{n_p}$ are the state, control and parameter variables, $\hat{x}_0$ the initial condition and $f$ an uniformly Lipschitz continuous function in $x,u$ and continuous in $t$. $h$ and $h_f$ are the optimization objectives, $r$ the path constraints and $l$ the boundary conditions.

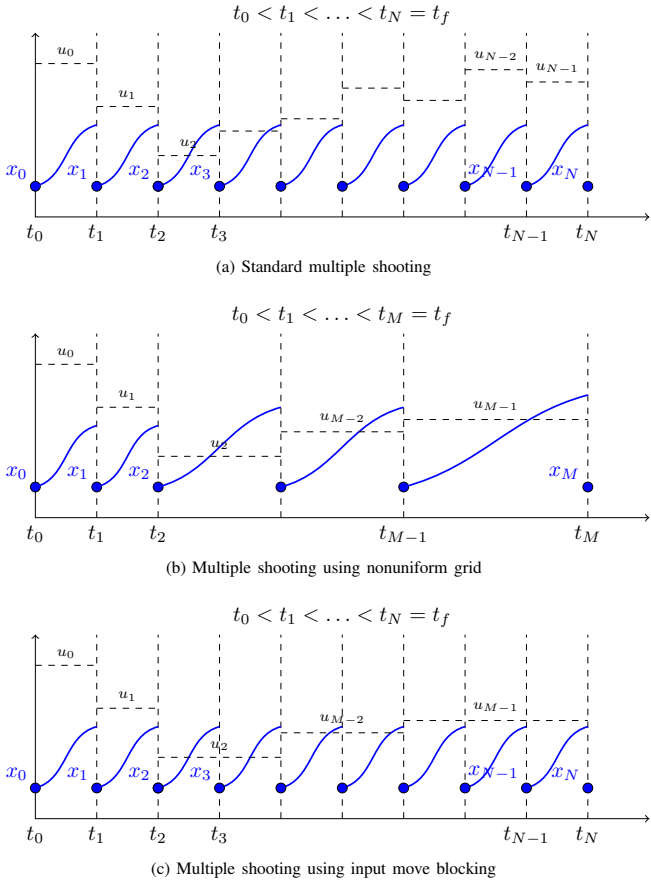
\begin{figure}[htb]
\centering
\begin{adjustbox}{max width=1\linewidth}
\subfloat[Standard multiple shooting]{
\begin{tikzpicture}
\draw [<->](0,3)--(0,0)--(10,0);
\draw [dashed, ultra thin] (1,0)--(1,3);
\node [below] at (0,0) {$t_0$};
\draw [dashed, ultra thin] (2,0)--(2,3);
\node [below] at (1,0) {$t_1$};
\draw [dashed, ultra thin] (3,0)--(3,3);
\node [below] at (2,0) {$t_2$};
\draw [dashed, ultra thin] (4,0)--(4,3);
\node [below] at (3,0) {$t_3$};
\draw [dashed, ultra thin] (5,0)--(5,3);
\draw [dashed, ultra thin] (6,0)--(6,3);
\draw [dashed, ultra thin] (7,0)--(7,3);
\draw [dashed, ultra thin] (8,0)--(8,3);
\draw [dashed, ultra thin] (9,0)--(9,3);
\node [below] at (8,0) {$t_{N-1}$};
\node [below] at (9,0) {$t_N$};
\draw [fill=blue] (0,0.5) circle [radius=0.08];
\draw [thick, blue] (0,0.5) to [out=15,in=195] (1,1.5);
\draw [fill=blue] (1,0.5) circle [radius=0.08];
\draw [thick, blue] (1,0.5) to [out=15,in=195] (2,1.5);
\draw [fill=blue] (2,0.5) circle [radius=0.08];
\draw [thick, blue] (2,0.5) to [out=15,in=195] (3,1.5);
\draw [fill=blue] (2,0.5) circle [radius=0.08];
\draw [thick, blue] (3,0.5) to [out=15,in=195] (4,1.5);
\draw [fill=blue] (3,0.5) circle [radius=0.08];
\draw [thick, blue] (4,0.5) to [out=15,in=195] (5,1.5);
\draw [fill=blue] (4,0.5) circle [radius=0.08];
\draw [thick, blue] (5,0.5) to [out=15,in=195] (6,1.5);
\draw [fill=blue] (5,0.5) circle [radius=0.08];
\draw [thick, blue] (6,0.5) to [out=15,in=195] (7,1.5);
\draw [fill=blue] (6,0.5) circle [radius=0.08];
\draw [thick, blue] (7,0.5) to [out=15,in=195] (8,1.5);
\draw [fill=blue] (7,0.5) circle [radius=0.08];
\draw [thick, blue] (8,0.5) to [out=15,in=195] (9,1.5);
\draw [fill=blue] (8,0.5) circle [radius=0.08];
\draw [fill=blue] (9,0.5) circle [radius=0.08];

\node [above] at (5,3.0) {$t_0<t_1<\ldots<t_N=t_f$};
\node [above left, blue] at (0,0.5) {$x_0$}; 
\node [above left, blue] at (1,0.5) {$x_1$}; 
\node [above left, blue] at (2,0.5) {$x_2$}; 
\node [above left, blue] at (3,0.5) {$x_3$};
\node [above left, blue] at (8,0.5) {$x_{N-1}$};
\node [above left, blue] at (9,0.5) {$x_{N}$};

\draw [dashed] (0,2.5)--(1,2.5);
\node [above] at (0.5,2.5) {\scriptsize $u_0$};
\draw [dashed] (1,1.8)--(2,1.8);
\node [above] at (1.5,1.8) {\scriptsize $u_1$};
\draw [dashed] (2,1.0)--(3,1.0);
\node [above] at (2.5,1.0) {\scriptsize $u_2$};
\draw [dashed] (3,1.4)--(4,1.4);
\draw [dashed] (4,1.6)--(5,1.6);
\draw [dashed] (5,2.1)--(6,2.1);
\draw [dashed] (6,1.9)--(7,1.9);
\draw [dashed] (7,2.4)--(8,2.4);
\node [above] at (7.5,2.4) {\scriptsize $u_{N-2}$};
\draw [dashed] (8,2.2)--(9,2.2);
\node [above] at (8.5,2.2) {\scriptsize $u_{N-1}$};
\end{tikzpicture}
\label{standardMPC}}
\end{adjustbox}

\vfill

\begin{adjustbox}{max width=1\linewidth}
\subfloat[Multiple shooting using nonuniform grid]{
\begin{tikzpicture}
\draw [<->](0,3)--(0,0)--(10,0);
\draw [dashed, ultra thin] (1,0)--(1,3);
\node [below] at (0,0) {$t_0$};
\draw [dashed, ultra thin] (2,0)--(2,3);
\node [below] at (1,0) {$t_1$};
\draw [dashed, ultra thin] (4,0)--(4,3);
\node [below] at (2,0) {$t_2$};
\draw [dashed, ultra thin] (6,0)--(6,3);
\draw [dashed, ultra thin] (9,0)--(9,3);
\node [below] at (6,0) {$t_{M-1}$};
\node [below] at (9,0) {$t_M$};
\draw [fill=blue] (0,0.5) circle [radius=0.08];
\draw [thick, blue] (0,0.5) to [out=15,in=195] (1,1.5);
\draw [fill=blue] (1,0.5) circle [radius=0.08];
\draw [thick, blue] (1,0.5) to [out=15,in=195] (2,1.5);
\draw [fill=blue] (2,0.5) circle [radius=0.08];
\draw [thick, blue] (2,0.5) to [out=15,in=195] (4,1.8);
\draw [fill=blue] (4,0.5) circle [radius=0.08];
\draw [thick, blue] (4,0.5) to [out=15,in=195] (6,1.8);
\draw [fill=blue] (6,0.5) circle [radius=0.08];
\draw [thick, blue] (6,0.5) to [out=15,in=195] (9,2.0);
\draw [fill=blue] (9,0.5) circle [radius=0.08];

\node [above] at (5,3.0) {$t_0<t_1<\ldots<t_M=t_f$};
\node [above left, blue] at (0,0.5) {$x_0$}; 
\node [above left, blue] at (1,0.5) {$x_1$}; 
\node [above left, blue] at (2,0.5) {$x_2$};
\node [above left, blue] at (9,0.5) {$x_{M}$};

\draw [dashed] (0,2.5)--(1,2.5);
\node [above] at (0.5,2.5) {\scriptsize $u_0$};
\draw [dashed] (1,1.8)--(2,1.8);
\node [above] at (1.5,1.8) {\scriptsize $u_1$};
\draw [dashed] (2,1.0)--(4,1.0);
\node [above] at (3,1.0) {\scriptsize $u_2$};
\draw [dashed] (4,1.4)--(6,1.4);
\draw [dashed] (6,1.6)--(9,1.6);
\node [above] at (5,1.4) {\scriptsize $u_{M-2}$};
\node [above] at (7.5,1.6) {\scriptsize $u_{M-1}$};
\end{tikzpicture}
\label{nonuniformgridMPC}}
\end{adjustbox}

\vfill

\begin{adjustbox}{max width=1\linewidth}
\subfloat[Multiple shooting using input move blocking]{
\begin{tikzpicture}
\draw [<->](0,3)--(0,0)--(10,0);
\draw [dashed, ultra thin] (1,0)--(1,3);
\node [below] at (0,0) {$t_0$};
\draw [dashed, ultra thin] (2,0)--(2,3);
\node [below] at (1,0) {$t_1$};
\draw [dashed, ultra thin] (3,0)--(3,3);
\node [below] at (2,0) {$t_2$};
\draw [dashed, ultra thin] (4,0)--(4,3);
\node [below] at (3,0) {$t_3$};
\draw [dashed, ultra thin] (5,0)--(5,3);
\draw [dashed, ultra thin] (6,0)--(6,3);
\draw [dashed, ultra thin] (7,0)--(7,3);
\draw [dashed, ultra thin] (8,0)--(8,3);
\draw [dashed, ultra thin] (9,0)--(9,3);
\node [below] at (8,0) {$t_{N-1}$};
\node [below] at (9,0) {$t_N$};
\draw [fill=blue] (0,0.5) circle [radius=0.08];
\draw [thick, blue] (0,0.5) to [out=15,in=195] (1,1.5);
\draw [fill=blue] (1,0.5) circle [radius=0.08];
\draw [thick, blue] (1,0.5) to [out=15,in=195] (2,1.5);
\draw [fill=blue] (2,0.5) circle [radius=0.08];
\draw [thick, blue] (2,0.5) to [out=15,in=195] (3,1.5);
\draw [fill=blue] (2,0.5) circle [radius=0.08];
\draw [thick, blue] (3,0.5) to [out=15,in=195] (4,1.5);
\draw [fill=blue] (3,0.5) circle [radius=0.08];
\draw [thick, blue] (4,0.5) to [out=15,in=195] (5,1.5);
\draw [fill=blue] (4,0.5) circle [radius=0.08];
\draw [thick, blue] (5,0.5) to [out=15,in=195] (6,1.5);
\draw [fill=blue] (5,0.5) circle [radius=0.08];
\draw [thick, blue] (6,0.5) to [out=15,in=195] (7,1.5);
\draw [fill=blue] (6,0.5) circle [radius=0.08];
\draw [thick, blue] (7,0.5) to [out=15,in=195] (8,1.5);
\draw [fill=blue] (7,0.5) circle [radius=0.08];
\draw [thick, blue] (8,0.5) to [out=15,in=195] (9,1.5);
\draw [fill=blue] (8,0.5) circle [radius=0.08];
\draw [fill=blue] (9,0.5) circle [radius=0.08];

\node [above] at (5,3.0) {$t_0<t_1<\ldots<t_N=t_f$};
\node [above left, blue] at (0,0.5) {$x_0$}; 
\node [above left, blue] at (1,0.5) {$x_1$}; 
\node [above left, blue] at (2,0.5) {$x_2$}; 
\node [above left, blue] at (3,0.5) {$x_3$};
\node [above left, blue] at (8,0.5) {$x_{N-1}$};
\node [above left, blue] at (9,0.5) {$x_{N}$};

\draw [dashed] (0,2.5)--(1,2.5);
\node [above] at (0.5,2.5) {\scriptsize $u_0$};
\draw [dashed] (1,1.8)--(2,1.8);
\node [above] at (1.5,1.8) {\scriptsize $u_1$};
\draw [dashed] (2,1.0)--(4,1.0);
\node [above] at (3,1.0) {\scriptsize $u_2$};
\draw [dashed] (4,1.4)--(6,1.4);
\draw [dashed] (6,1.6)--(9,1.6);
\node [above] at (5,1.4) {\scriptsize $u_{M-2}$};
\node [above] at (7.5,1.6) {\scriptsize $u_{M-1}$};
\end{tikzpicture}
\label{inputMBMPC}}
\end{adjustbox}

\caption{Comparison of three different parameterization strategies: (a) the state and control trajectories are parameterized into $N$ shooting intervals; (b) the state and control trajectories are parameterized into $M$ shooting intervals. The length of each interval is non-equidistant; (c) the state and control trajectories are parameterized into $N$ and $M$ shooting intervals, respectively.}
\label{comparison}
\end{figure}

The solution of \eqref{OCP} can be found by applying direct methods that employ finite dimensional parameterization to the OCP obtaining a NLP problem. Popular parameterization methods include direct single shooting, multiple shooting, and collocation. Pros and cons of different methods are discussed in \cite{binder2001introduction}. In particular, multiple shooting \cite{bock1984multiple} has proven its effectiveness in several real-time applications. In multiple shooting, the prediction horizon is usually divided into $N$ equidistant shooting intervals using time grid $\left [ t_k, t_{k+1} \right ]$, $k=0,\dots,N-1$ and $t_N=t_f$. The state trajectory is discretized at the $N+1$ discretization nodes and the control input is assumed piece-wise constant on each shooting interval (see Fig \ref{standardMPC}). By also parameterizing  the objective function and constraints, a NLP problem is formulated as
\begin{equation}\label{NLP}
\begin{aligned}
\min_{\mathbf{x},\mathbf{u}} \quad& \sum_{k=0}^{N-1} h(t_k,x_k,u_k;p)+h_N(x_N;p)\\
s.t.\quad & 0 = \hat{x}_0- x_0\\
& 0 =\phi(t_k,x_k,u_k;p)-x_{k+1},\,\forall  k=0,1,\dots, N-1,\\
& 0 \geq r(x_k,u_k;p), \, \forall  k=0,1,\dots, N-1,\\
& 0 \geq l(x_N;p).
\end{aligned}
\end{equation}
where $\phi$ is a numerical integration operator that solves the following initial value problem (IVP) and return the solution at $t_{k+1}$.
\begin{equation}\label{IVP}
0=f(\dot{x}(t), x(t),u(t),t),\quad x(0)=x_k.
\end{equation}
We define 
\begin{equation}
    \begin{aligned}
    &\mathbf{x}= \left [ x_0^\top, x_1^\top,\dots, x_N^\top\right ]^\top,\\
    &\mathbf{u}= \left [ u_0^\top, u_1^\top,\dots, u_{N-1}^\top\right ]^\top,
    \end{aligned}
\end{equation}
as the discrete state and control variables. Problem \eqref{NLP} can be solved by means of SQP and the resulting QP problem, obtained by linearizing \eqref{NLP} at $(\mathbf{x}^i, \mathbf{u}^i)$, reads:
\begin{equation}\label{QP}
\begin{aligned}
\min_{\Delta \mathbf{x},\Delta \mathbf{u}} \quad& \sum_{k=0}^{N-1} \left (\frac{1}{2} 
\begin{bmatrix}
\Delta x_k \\
\Delta u_k
\end{bmatrix}^\top H_k^i
\begin{bmatrix}
\Delta x_k \\
\Delta u_k
\end{bmatrix}+g_k^{i^\top}
\begin{bmatrix}
\Delta x_k \\
\Delta u_k
\end{bmatrix} \right ) \\
&+\frac{1}{2} \Delta x_N^\top H_N^i \Delta x_N^\top+g_N^{i^\top} \Delta x_N\\
s.t.\quad & \Delta x_0= \hat{x}_0-x_0^i\\
& \Delta x_{k+1}=A_{k}^i \Delta x_{k}+ B_{k}^i \Delta u_{k} + d_{k}^i,  \\
& C_k^i
\begin{bmatrix}
\Delta x_k \\
\Delta u_k
\end{bmatrix} + c_k^i\leq 0, \quad \forall  k=0,1,\dots, N-1,\\
& C_N^i \Delta x_N + c_N^i \leq 0,
\end{aligned}
\end{equation}
where $\Delta \mathbf{x}=\mathbf{x}-\mathbf{x}^i, \Delta \mathbf{u}=\mathbf{u}-\mathbf{u}^i$. The matrix $H_k^i$ is the k-th block of the Hessian associated to the Lagrangian of \eqref{NLP} and $g_k^i$ is the k-th sub-vector of the gradient of the objective function. Matrices $A_k^i=\frac{\partial \phi}{\partial x_k}(x_k^i,u_k^i), \,B_k^i=\frac{\partial \phi}{\partial u_k}(x_k^i,u_k^i)$ are sensitivities of dynamics at shooting node $k$ and the matrix $C_k^i=\frac{\partial r}{\partial (x_k, u_k)}(x_k^i,u_k^i)$ is the $k$-th block of the Jacobian matrix of the constraint. 

Problem \eqref{QP} has a particular multi-stage structure with $(N+1)n_x+Nn_u$ decision variables and can be solved by structure exploiting or sparse solvers \cite{zanelli2017forces, stellato2018osqp}. Such solvers are tailored for multi-stage problems with a complexity growing linearly with the number of discretization nodes $N$. Note that the state variables depend on control variables hence problem \eqref{QP} has $Nn_u$ DoFs. However, interior point methods cannot straightforwardly exploit warm-start strategies which are often desired for real-time NMPC applications. An alternative for solving problem \eqref{QP} is to exploit the multi-stage structure by eliminating the state variables from the equality constraint, leading to a condensed QP problem with $N n_u$ decision variables:
\begin{equation}\label{QP_full_condensed}
\begin{aligned}
\min_{\Delta \mathbf{u}} \quad& \frac{1}{2} {\Delta \mathbf{u}}^\top H_c \Delta \mathbf{u}+ g_c^\top \Delta \mathbf{u}\\
s.t.\quad & C_c \Delta \mathbf{u} + c_c \leq 0,
\end{aligned}
\end{equation}
The condensed problem \eqref{QP_full_condensed} can be efficiently solved by active-set methods using warm-start. Efficient algorithms for the condensing step, i.e. the computation of $H_c, g_c$ and so on, have been proposed in \cite{frison2013fast, andersson2013general}. The computation complexity of these condensing algorithms are $\mathcal{O}(N^2)$, mostly dominated by the computation of $H_c$ with $\mathcal{O}(N^2 n_x^2n_u+N^2 n_xn_u^2)$ floating point operations (FLOPs). 

In this paper, we adopt the RTI scheme which solves only one QP problem at each sampling instant \cite{diehl2002real}. The solution of \eqref{QP} or \eqref{QP_full_condensed}, is updated using a single, full Newton step without achieving local optimality of \eqref{NLP}. Despite its suboptimality, RTI scheme has been widely applied for real-time applications with fast changing dynamics. 

\section{Move Blocking strategy}\label{sec3}

Since the DoFs of \eqref{QP_full_condensed} is linear in $N$ and the complexity of condensing algorithm is quadratic in $N$, computation burden rises quickly when using a long prediction horizon, i.e. a large $N$. Therefore, for real-time applications, it is often recommended to employ condensing and dense QP solvers when $N$ is small. A way to reduce this complexity is to use input MB strategy, which fixes the control inputs to be constant over a certain number of successive time intervals over the prediction horizon. This can be achieved by adding equality constraints to \eqref{QP} to fix the values of elements in the vector of controls $\mathbf{u}$ \cite{cagienard2007move}. A typical way is to define $\Delta \hat{\mathbf{u}}$ by 
\begin{equation} 
\begin{aligned}
\Delta \mathbf{u}&=T \Delta\hat{\mathbf{u}}=(T_b \otimes I_{n_u \times n_u}) \Delta\hat{\mathbf{u}}\\
&= 
\begin{bmatrix}
E_0& &&\\
&E_1&&\\
&&\ddots&\\
&&&E_{M-1}
\end{bmatrix}
\begin{bmatrix}
\Delta\hat{u}_0\\
\Delta\hat{u}_1\\
\vdots\\
\Delta\hat{u}_{M-1}
\end{bmatrix},
\end{aligned}
\label{MB constraint}
\end{equation}
where $E_j$, for $j=0,...,M-1$ consist of $N_j$ vertically stacked identity matrices of size $n_u$, and $\Delta\hat{\mathbf{u}} = [\Delta\hat{u}_0^\top, \Delta\hat{u}_1^\top,...\Delta\hat{u}_{M-1}^\top]^\top$ with $M<N$ represents the sequence of the new control input, each applied for $N_j$ successive shooting intervals.
The number of columns of matrix $T$ is the number of DoFs (i.e. $M$) for the control sequence and the number of rows of $E_j$ is the length of each input block (i.e. $N_j$) over the prediction horizon. Also define $I = [I_0,\ldots,I_M]$ as the vector consisting of the starting index of each input block, which can be computed by
\begin{equation}
\begin{aligned}
    &I_0=0, \\
    &I_j=\sum_{k=0}^{j-1}N_k, \quad j=1,\ldots,M-1,\\
    &I_M=\sum_{k=0}^{M-1}N_k = N.
\end{aligned}
\label{index}
\end{equation}

\subsection{Embed MB into multiple shooting}\label{proposed algorithm}
Given problem (\ref{QP_full_condensed}), an intuitive way to introduce MB strategy is following the linear MPC MB scheme by adding the equality constraint (\ref{MB constraint}) to the condensed problem \eqref{QP_full_condensed}, obtaining the following QP problem:
\begin{equation}
\begin{aligned}
\min_{\Delta{\hat{\mathbf{u}}}} \quad& \frac{1}{2} {\Delta \hat{\mathbf{u}}}^\top \hat{H}_c {\Delta\hat{\mathbf{u}}}+ \hat{g}_c^\top \Delta \hat{\mathbf{u}}\\
s.t.\quad & \hat{C}_c \Delta \hat{\mathbf{u}} + c_c \leq 0,\\
\end{aligned}
\label{QP_MB_1}
\end{equation}
where $\hat{H}_c=T^\top H_c T$, $\hat{g}_c= T^\top g_c$ and  $\hat{C}_c=C_c T$. Here the computation of dense matrices such as $\hat{H}_c$, which involves the multiplication by $T$, can be performed in an efficient way by exploiting the particular structure of $T$. However, such computations require a condensing step with a complexity of $\mathcal{O}(N^2)$ before MB is introduced. Therefore, MB would be an additional computation burden to the standard NMPC scheme.

\begin{figure}[!ht]
\centering
\begin{tikzpicture}[node distance=1.5cm]
\node (st1) [block] {multiple shooting};
\node (mb1) [block, right of=st1, xshift=3cm, fill=yellow] {multiple shooting with MB};
\node (st2) [block, below of=st1] {condensing $\mathcal{O}(N^2)$};
\node (mb2) [block, right of=st2, xshift=3cm, fill=yellow] {tailored condensing $\mathcal{O}(N)$};
\node (st3) [block, below of=st2] {adding MB constraint \eqref{MB constraint}};

\node (st4) [block, below of=st3] {solving QP (DoFs=$Mn_u$)};
\node (mb3) [block, right of=st4, xshift=3cm, fill=yellow] {solving QP (DoFs=$Mn_u$)};

\draw [arrow] (st1) -- (st2);
\draw [arrow] (mb1) -- (mb2);
\draw [arrow] (st2) -- (st3);
\draw [arrow] (mb2) -- (mb3);
\draw [arrow] (st3) -- (st4);

\end{tikzpicture}
\caption{An illustration of the two options to introduce MB into multiple shooting based nonlinear MPC. Red: the intuitive way; Yellow: the proposed MB scheme}
\label{fig:MBscheme}
\end{figure}
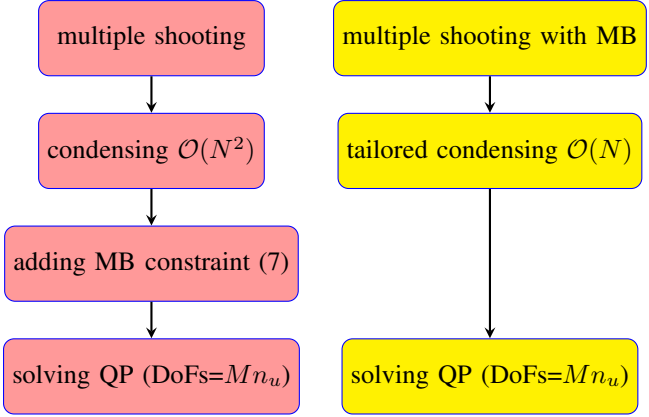

A more efficient alternative is to integrate the MB strategy during the multiple shooting phase. To do so, we initialize \eqref{NLP} using $\hat{\mathbf{u}}$, where the same input is applied to solve the differential equation \eqref{IVP} over several consecutive intervals (see Fig. \ref{inputMBMPC}). As a result, we obtain a QP problem with exactly the same structure as \eqref{QP} but has less DoFs: 
\begin{equation}\label{QP_MB}
\begin{aligned}
\min_{\Delta{\mathbf{x}},\Delta{\mathbf{x}}} \quad& \sum_{k=0}^{N-1} \left (\frac{1}{2} 
\begin{bmatrix}
\Delta x_k \\
\Delta \hat{u}_j
\end{bmatrix}^T H_k
\begin{bmatrix}
\Delta x_k \\
\Delta \hat{u}_j
\end{bmatrix}+{g_k}^T
\begin{bmatrix}
\Delta x_k \\
\Delta \hat{u}_j
\end{bmatrix} \right ) \\
&+\frac{1}{2} \Delta x_N^T H_N \Delta x_N^T+{g_N}^T \Delta x_N\\
s.t.\quad & \Delta x_0= \hat{x}_0-x_0\\
& \Delta x_{k+1}=A_{k} \Delta x_{k}+ B_{k} \Delta \hat{u}_{j} + d_{k},  \\
& C_k
\begin{bmatrix}
\Delta x_k \\
\Delta \hat{u}_j
\end{bmatrix} + c_k\leq 0, \quad \forall  k=0,1,\dots, N-1,\\
& C_N \Delta x_N + c_N \leq 0,
\end{aligned}
\end{equation}
where $\Delta \hat{u}_j$ is applied when $k\in[I_{j}, I_{j+1}), \forall j=0,1,\dots, M-1$, $A_k= \frac{\partial \phi}{\partial x}(x_k,\hat{u}_j)$, $B_k= \frac{\partial \phi}{\partial u}(x_k,\hat{u}_j)$, $d_k= \phi (t_k,x_k,\hat{u}_j;p)-x_{k+1}$. Since in \eqref{QP_MB} there are still $N+1$ stages in the cost function and $N$ stages in constraints, the level of sparsity is maintained. Hence, no additional computation burden is introduced for embedding input MB into multiple shooting based NMPC algorithm. An illustration of the two options to introduce MB into nonlinear MPC is given in Fig. \ref{fig:MBscheme}.

\subsection{Comparison to nonuniform grid schemes}
Nonuniform grid schemes are closely related to MB and can also reduce the dimension of \eqref{NLP} by using a denser grid at part of the prediction horizon and a coarser one in the other part of the horizon \cite{quirynen2015multiple, yu2016stable, potena2018non}. For example, as shown in Fig. \ref{nonuniformgridMPC}, if $M$ non-equidistant shooting intervals are defined, the state and control vectors shrink to 
\begin{equation}
    \begin{aligned}
    &\mathbf{x}'= \left [ x_0^\top, x_1^\top,\dots, x_M^\top\right ]^\top,\\
    &\mathbf{u}'= \left [ u_0^\top, u_1^\top,\dots, u_{M-1}^\top\right ]^\top.
    \end{aligned}
\end{equation}
The resulting QP has the same structure as \eqref{QP} with $M+1$ stages. In fact, non-uniform grid schemes and input MB schemes belong to the same kind of parameterization limiting algorithms, where the former limits both state and input parameterizations and the latter only limits the input. As a result, MB has several advantages over nonuniform grid schemes, which are summarized in the following.
\begin{enumerate}
    \item A more accurate state trajectory can be obtained by using MB since the grid for state discretization is as precise as the uniform grid MPC. Nonuniform grid schemes suffer from a more coarse state discretization and less accurate predicting trajectory. Although such inaccuracy can be alleviated by using time-mesh refinement techniques \cite{paiva2015adaptive, potena2018non, lee2018mesh, paiva2018sufficient}, additional algorithm complexity and computation burden is introduced. The computational time is also time-varying due to different number of discretizing nodes at each sampling step.
    \item Path and state constraint fulfillment can be maintained using MB. Since in \eqref{NLP} constraints are only fulfilled exactly at the shooting nodes but not between them, nonuniform grid schemes would have a larger possibility for constraint violation due to less state discretization nodes.
    \item The proposed MB strategy has $N+1$ stages in cost function as in \eqref{NLP}, while nonuniform grid schemes has only $M$ stages. As a result, nonuniform grid schemes usually need modified weights for each stage and such modification is not trivial to obtain the same optimization performance as that of uniform grid NMPC algorithms.
    \item For $M \ll N$, the level of sparsity of \eqref{NLP} is preserved using the proposed MB while nonuniform grid schemes lose the sparsity. Such level of sparsity is considerably important to improve the convergence property and solution accuracy for solving \eqref{NLP} numerically \cite{quirynen2015multiple}. 
\end{enumerate}

\subsection{Tailored condensing}
We propose a tailored condensing algorithm that exploits the sparsity structure of \eqref{QP_MB} to obtain \eqref{QP_MB_1}. By exploiting the reduced DoFs of \eqref{QP_MB}, the computational complexity of the condensing step is reduced from $\mathcal{O}(N^2)$ to $\mathcal{O}(NM)$. 

\begin{algorithm}[t] 
\caption{Calculation of  $\hat{G}=GT$ with complexity $\mathcal{O}(N M)$}
\label{G}
\begin{algorithmic}
    \Initialize{$\hat{G} \gets 0_{Nn_x \times Mn_u}$}
    \For{$i = 0,..., M-1$}
	 \State $\hat{G}[I_i,i] \gets B_{I_i} $
        \For{$j = I_i+1,..., N-1$}
	      \If{$j < I_{i+1}$}
       	\State $\hat{G}[j,i] \gets A_j \quad \hat{G}[j-1,i]+B_j$
    \Else
        \State $\hat{G}[j,i] \gets A_j \quad \hat{G}[j-1,i]$
    \EndIf
    \EndFor
    \EndFor
\end{algorithmic}
\end{algorithm}

\begin{algorithm}[ht]
\caption{Calculation of $\hat{H}_c=T^\top H_c T$ with complexity $\mathcal{O}(N M)$} 
\label{H_c} 
\begin{algorithmic} 
    \Initialize{$\hat{H}_c \gets 0_{Mn_u\times Mn_u}, \, \hat{H}_{tmp} \gets 0_{Nn_u\times Mn_u}, \,R_{tmp}\gets 0_{n_u \times n_u}$}
    \For{$i = 0,...,M-1$}
	    \State $ W_N \gets Q_N \hat{G}[N-1,i]$
    	\For{$k = N-1,...,I_i+1$}
	            \State ${H}_{tmp}[k,i] \gets S_k^\top\hat{G}[k-1,i] + B_k^\top  W_{k+1}$
	    \State ${W}_k \gets Q_k \hat{G}[k-1,i] + A_k^\top W_{k+1}$
	    \EndFor
	    \State ${H}_{tmp}[I_i,i]  \gets B_{I_i}^\top W_{I_i+1}$
    \EndFor
    
    \State $k \gets 0$
    
    \For{$i$ = $0,...,N-1$}
        \State $\hat{H}_c[k,0:M-1] += H_{tmp} [i,0:M-1]$
        \State $R_{tmp} +=R_i$
        \If{$i+1=I_{k+1}$}
            \State $\hat{H}_c[k,k] += R_{tmp}$
            \State $k++$
            \State $R_{tmp}\gets 0$
        \EndIf
    \EndFor
\end{algorithmic}
\end{algorithm}

According to the equality constraint in \eqref{QP} and \eqref{QP_MB}, the state variable $\Delta \mathbf{x}$ can be expressed by the input variable $\Delta \mathbf{u}$ and $\Delta \mathbf{\hat{u}}$. Let us define matrices $G,\hat{G}$ that satisfy
\begin{equation}
    \Delta \mathbf{x} = G \Delta \mathbf{u}+L = \hat{G}\Delta \hat{\mathbf{u}}+L,
\end{equation}
where 
\begin{equation}
    G=\begin{bmatrix}
    G_{0,0} \\
    G_{1,0} & G_{1,1} \\
    \vdots & & \ddots \\
    G_{N-1,0} & G_{N-1,1} & \cdots & G_{N-1,N-1},
    \end{bmatrix},
\end{equation}
with $N$ row and column blocks and $L$ the corresponding residual vector, which can be computed as in \cite{andersson2013general}. Similarly, we have
\begin{equation}
    \hat{G}=\begin{bmatrix}
    \hat{G}_{0,0} \\
    \hat{G}_{1,0} & \hat{G}_{1,1} \\
    \vdots & & \ddots \\
    \hat{G}_{N-1,0} & \hat{G}_{N-1,1} & \cdots & \hat{G}_{N-1,M-1},
    \end{bmatrix},
\end{equation}
has $N$ row and $M$ column blocks. In \eqref{QP_MB}, $\hat{H}_k$ has the form of 
\begin{equation}\label{Hessian detail}
\hat{H}_k=\begin{bmatrix}
Q_k & S_k \\
S_k^\top & R_k
\end{bmatrix}.
\end{equation}
Inspired by \cite{andersson2013general}, we propose modified condensing algorithms that exploit the reduced DoF in \eqref{QP_MB} without the need to explicitly build the matrix $T$. Algorithm \ref{G} and Algorithm \ref{H_c} present the computation of $\hat{G}$ and $\hat{H}_c$ which are leading factors in condensing algorithms. In these algorithms, $X[i,j]$ is the block matrix at the $i^{th}$ row and $j^{th}$ column of $X$. The computational complexity for computing $\hat{H}_c$ is $\mathcal{O}(NMn_x^2n_u+NMn_xn_u^2)$ FLOPs.

It is worth noticing that, since only the first control component $u^*_0$ is fed to the system, a coarse discretization of intervals far from the current time instant in the prediction horizon may achieve the same degree of solution accuracy \cite{paiva2018sufficient}. In addition, as indicated in the linear MPC case \cite{shekhar2015optimal}, there exists an optimal block structure that maintains the size of ROA of MPC. Hence, the number of DoFs is not monotonically related to the solution accuracy nor the size of ROA, which is demonstrated by simulations in Section \ref{sec4}. 

\begin{figure*}[!htb]
    \centering
    \includegraphics[width=.8\textwidth]{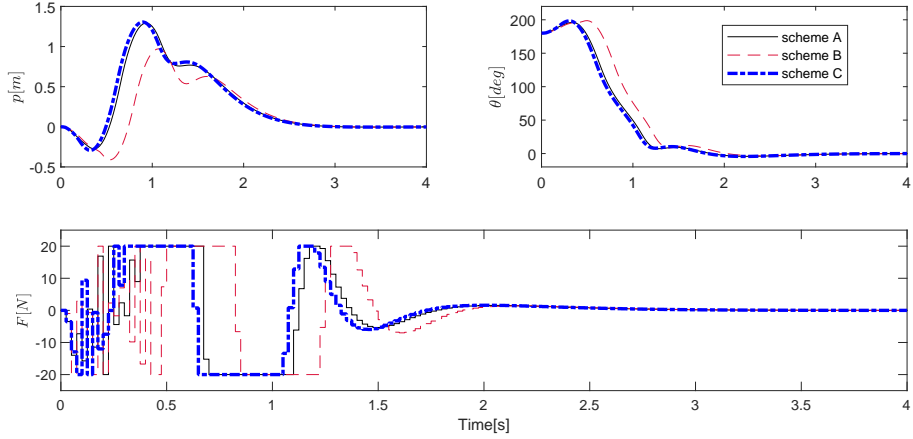}
    \caption{Closed-loop state and control trajectories using the three NMPC schemes for problem \eqref{invertedpendulumMPC}}
    \label{traj_compare}
\end{figure*}

\subsection{Convergence and Stability}
The off-line convergence property of the proposed input MB NMPC scheme can be analyzed in the framework of Newton-type methods. For simplicity, let us consider problem \eqref{NLP} without inequality constraints. Without loss of generality, we choose a block structure with $M=2$, a free $u_0$ and $u_1=u_2=\ldots=u_{N-1}$. We assume such a structure is recursively feasible, convergent and stable for a specific system.

Re-write \eqref{NLP} in the following compact form:
\begin{equation}\label{NLP1}
    \begin{aligned}
        \min_{x,u,\lambda} \quad& F(x,u)\\
        s.t. \quad & V(x,u)=0,\\
             &B^\top u=0,
    \end{aligned}
\end{equation}
where 
\begin{equation}
    B^\top u=\begin{bmatrix}
    u_1-u_2\\
    u_2-u_3\\
    \vdots\\
    u_{N-2}-u_{N-1}
    \end{bmatrix}.
\end{equation}
The Lagrangian of \eqref{NLP1} is
\begin{equation}
    \tilde{\mathcal{L}}(x,u,\lambda,\lambda_M) = \mathcal{L}(x,u,\lambda)+\lambda_M^\top B^\top u
\end{equation}
The necessary condition of optimality is:
\begin{equation}
    \nabla\tilde{\mathcal{L}}= \begin{bmatrix}
    \nabla_u \mathcal{L}+B\lambda_M\\
    \nabla_x \mathcal{L}\\
    V\\
    B^\top u
    \end{bmatrix}=0.
\end{equation}
Assume there exist matrices $Q_1,Q_2$ such that
\begin{equation}
    Q_1^\top B=I, \quad Q_2^\top B=0.
\end{equation}
As a result, by multiplying $(Q1|Q2)^\top$ to the first element of $\nabla\tilde{\mathcal{L}}$ we get
\begin{equation}\label{MB Lang grad}
    \begin{bmatrix}
    Q_1^\top \nabla_u \mathcal{L}+\lambda_M\\
    Q_2^\top \nabla_u \mathcal{L}\\
    \nabla_x \mathcal{L}\\
    V\\
    B^\top u
    \end{bmatrix}=0.
\end{equation}
The additional multiplier $\lambda_M$ can be chosen to be
\begin{equation}
    \lambda_M = -Q_1^\top \nabla_u \mathcal{L},
\end{equation}
making the first component of \eqref{MB Lang grad} always zero. As a result, we can apply Theorem 5.4 in \cite{diehl2001real} to show that if the original unblocked NMPC is convergent, the input MB NMPC is also convergent with the same convergence rate.

In this paper, we have set the problem in the framework of RTI scheme \cite{diehl2002real} to achieve an efficient implementation. Since only one SQP iteration is performed at each sampling instant in the RTI scheme, the on-line convergence instead of off-line convergence has to be considered \cite{diehl2005nominal}. As shown in \cite{diehl2001real}, the RTI scheme can be considered as a standard NMPC algorithm with the first control element fixed to the value from the previous sampling instant. Therefore, the same result apply to the RTI case, where the convergence rate of input MB RTI is the same with the standard RTI.

For stability analysis, input MB NMPC can be considered as a sub-optimal NMPC algorithm which benefits from (sub-optimal) NMPC stability theories \cite{mayne2000constrained, pannocchia2011conditions, graichen2010stability, allan2017inherent}. 
As shown in \cite{paiva2018sufficient} (Theorem 5.1) and \cite{yu2016stable}, a stabilizing non-equidistantly discretized NMPC can be obtained by choosing a sufficiently long prediction horizon, properly shifting the optimal input trajectory, properly designing the cost function, the terminal cost and constraints. Since input MB is a subclass of non-equidistant discretization, the developed MB NMPC is closed-loop stable once the initial block structure satisfies the stability conditions.

It is worth noting that the block structure is a key design factor affecting the size of ROA in the presence of terminal constraints, as well as the matrices $Q_1,Q_2$ in \eqref{MB Lang grad}. Analyses aforementioned assume a block structure that satisfy the convergence and stability conditions, but do not point out how to find such one. This issue is handled in linear MPC by solving two mixed-integer linear programs (MILP) off-line \cite{shekhar2015optimal}, one for finding the minimal number of blocks to maintain the size of ROA and the other for finding the maximal size of ROA when the number of blocks is fixed. However, for NMPC, one need to solve mixed-integer nonlinear programs (MINLP) which is NP-hard and may need on-line computations, implying the optimal block structure being time-varying. In this paper we focus on computationally efficient algorithms, hence a block structure that leads to convergence and closed-loop stability is always assumed.

\section{Implementation and Numerical Examples}\label{sec4}
The proposed input MB scheme is implemented in MATMPC \cite{chen2018matmpc}, an open source software built in MATLAB for NMPC applications. MATMPC has a number of algorithmic modules and provides state-of-the-art computation performance while making the prototyping easy with limited programming knowledge. This is achieved by writing each module directly in MATLAB API for C. As a result, MATMPC modules can be compiled into MEX functions with performance comparable to plain C/C++ solvers. MATMPC has been successfully used in operating systems including WINDOWS, LINUX and OS X \cite{chen2017fast,chen2017inexact,chen2018adaptive,chen2018efficient}. 

In this section, two numerical examples are used to compare the numerical, control and computation performance of the three NMPC schemes in Fig. \ref{comparison}, namely 1) standard NMPC (refer to scheme A); 2) NMPC with nonuniform grid discretization (refer to scheme B); 3) NMPC with the proposed input MB (refer to scheme C). The first example is a classical inverted pendulum system and the second is a nontrivial dynamic driving simulator. The simulation is performed in MATLAB in WINDOWS 10 using MATMPC \cite{chen2018matmpc} as the on-line solver, on a PC with Intel core i7-4790 running at 3.60GHz. The QP problem is solved by qpOASES \cite{ferreau2014qpoases}.

\begin{figure*}[!htb]
    \centering
    \includegraphics[width=.8\textwidth]{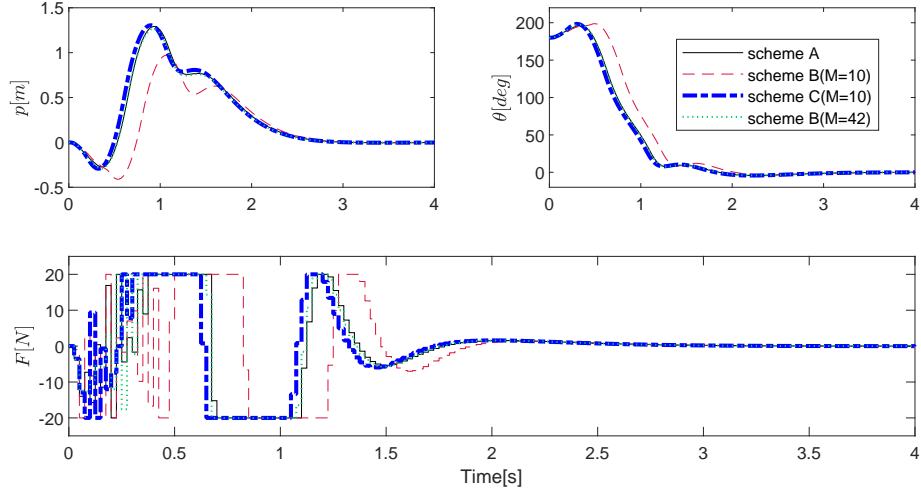}
    \caption{Closed-loop state and control trajectories using the three NMPC schemes for problem \eqref{invertedpendulumMPC}. Two different number of discretization nodes for scheme B are used.}
    \label{traj_compare_new}
\end{figure*}

\subsection{Control of an inverted pendulum}
\begin{figure}[!ht]
    \centering
    \includegraphics[width=.5\textwidth]{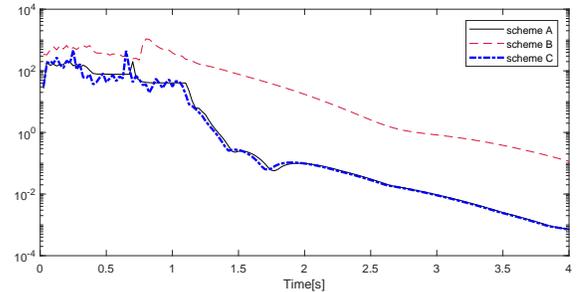}
    \caption{Closed-loop KKT values using the three NMPC schemes for problem \eqref{invertedpendulumMPC}.}
    \label{kkt_compare}
\end{figure}

An inverted pendulum is mounted on top of a cart and can roll up to 360 degrees. The dynamic model is given by
\begin{equation}
\begin{aligned}\label{inverted pendulum}
\ddot{p}&=\frac{-m_1l\sin(\theta)\dot{\theta}^2+m_1g\cos(\theta)\sin(\theta)+u}{m_2+m_1-m_1(\cos(\theta))^2},\\
\ddot{\theta}&=\frac{1}{l(m_2+m_1-m_1(\cos(\theta))^2)}(u\cos(\theta)\\
&-m_1l\cos(\theta)\sin(\theta)\dot{\theta}^2\\
&+(m_2+m_1)g\sin(\theta)),
\end{aligned}
\end{equation}
where $p,\theta$ are the cart position and swinging angle, respectively, and $u$ is the control force acting on the cart. The NMPC problem is formulated as
\begin{equation}
    \begin{aligned}\label{invertedpendulumMPC}
        \min_{x(\cdot),u(\cdot)} \quad &\int_{t_0}^{t_f} (\norm{x(t)}_{Q}^2+\norm{u(t)}_{R}^2)\text{d}t + \norm{x(t_f)}_{Q_f}^2 \\
        s.t. \quad &x(t_0) = \hat{x}_0, \\
            & \text{dynamics given in } (9) \\
        & -2\leq p(t)\leq 2,\\
        & -20\leq u(t)\leq 20
    \end{aligned}
\end{equation}
where $x(t)=[p(t),\theta(t),\dot{p}(t),\dot{\theta}(t)]^\top$. The model and values of parameters $m_1,m_2,l,g$ are taken from \cite{quirynen2015autogenerating}. The control task is to drive the pendulum from the downward position ($x=[0,\pi,0,0]$) to the upward position ($x=[0,0,0,0]$). We evaluate the algorithm performance in the following configurations:

\begin{enumerate}
    \item \textit{Standard RTI (scheme A)}: with sampling and shooting interval time $T_s=25$ms and $N=80$. 
    \item \textit{Nonuniform grid (scheme B)}: the index of nonuniform grid is given by $I = [0, 1, 3, 6, 10, 15, 20, 35, 50, 65, 80]$. Hence the number of stages is $M=10$. The weights for the cost function are based on that used by scheme A, scaled by the length of each shooting interval.
    \item \textit{MB (scheme C)}: the index of input block is given by the same I as in scheme B. Hence, the number of stages is $N=80$ but the DoFs is $M=10$. The weights for the cost function are identical to that used by scheme A.
\end{enumerate}
Fig. \ref{traj_compare} shows the closed-loop trajectories using the three MPC schemes for the inverted pendulum example. The input MB scheme has a very similar control performance comparing to the standard MPC, without the need to change the cost function weights. The trajectory of nonuniform grid scheme deviate from that of the standard MPC and the control task is accomplished in a longer time. Its performance can be improved by adjusting the scaling factor of weights. However, the scaling of weights requires practical experience and is not trivial for many applications. Fig. \ref{kkt_compare} shows the closed-loop Karush–Kuhn–Tucker (KKT) value that reflects the degree of optimality of the NLP solution at each sampling time. Note that we implement all algorithms in RTI framework hence the KKT value is not affected by the number of iterations (since there is only one iteration), but the features of problem \eqref{QP}, e.g. level of sparsity. As sparsity structure is maintained, the proposed MB scheme has a KKT at the same level as the standard MPC using only one eighth of the DoFs. The nonuniform grid scheme has a much higher KKT due to sparsity lost, leading to more inaccurate solutions and poorer control performance.

\begin{figure}[!ht]
    \centering
    \includegraphics[width=.5\textwidth]{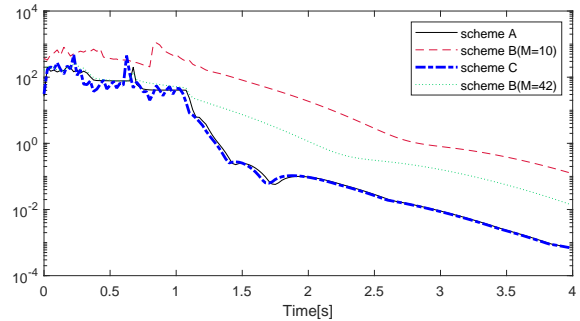}
    \caption{Closed-loop KKT values using the three NMPC schemes for problem \eqref{invertedpendulumMPC}. Two different number of discretization nodes for scheme B are used.}
    \label{kkt_compare_new}
\end{figure}

\begin{figure*}[!ht]
    \centering
    \subfloat[]{\includegraphics[width=.48\textwidth]{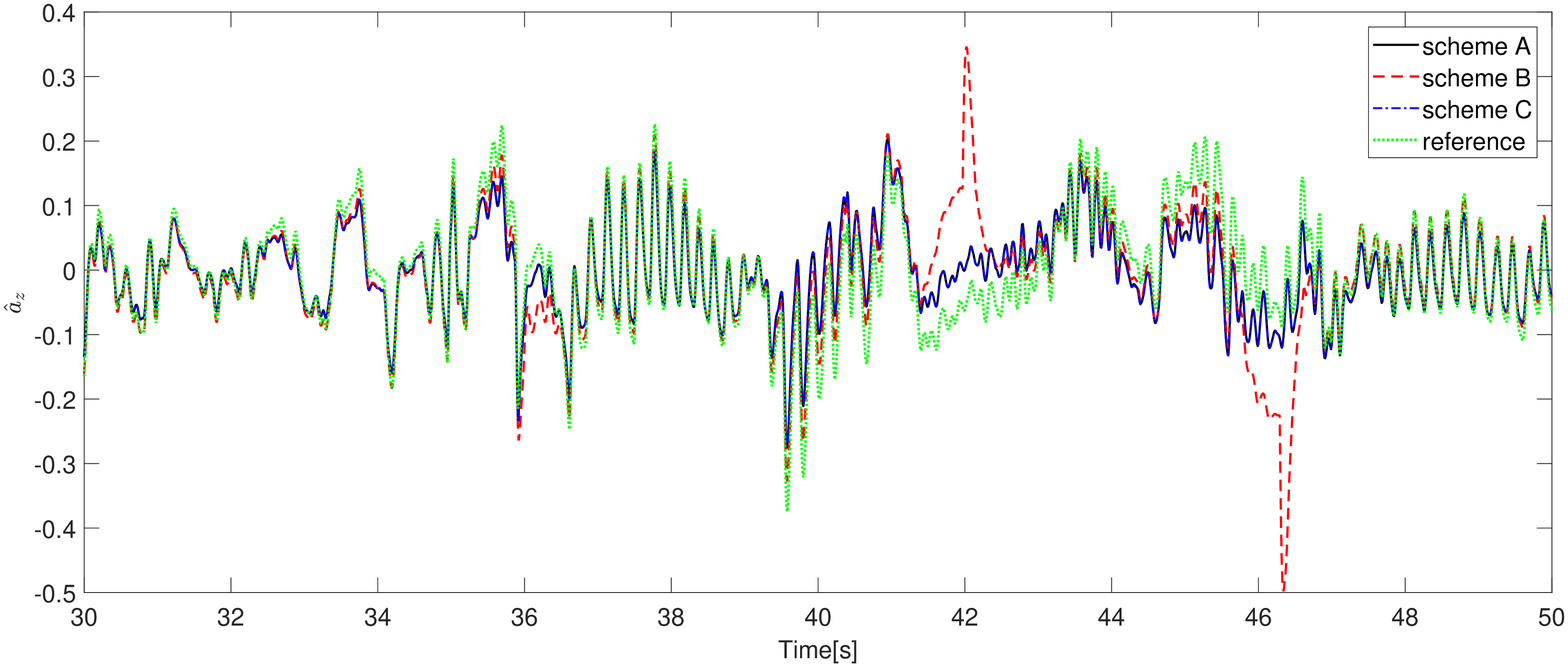}\label{traj_compare_DiM}}
    \quad
    \subfloat[]{\includegraphics[width=.48\textwidth]{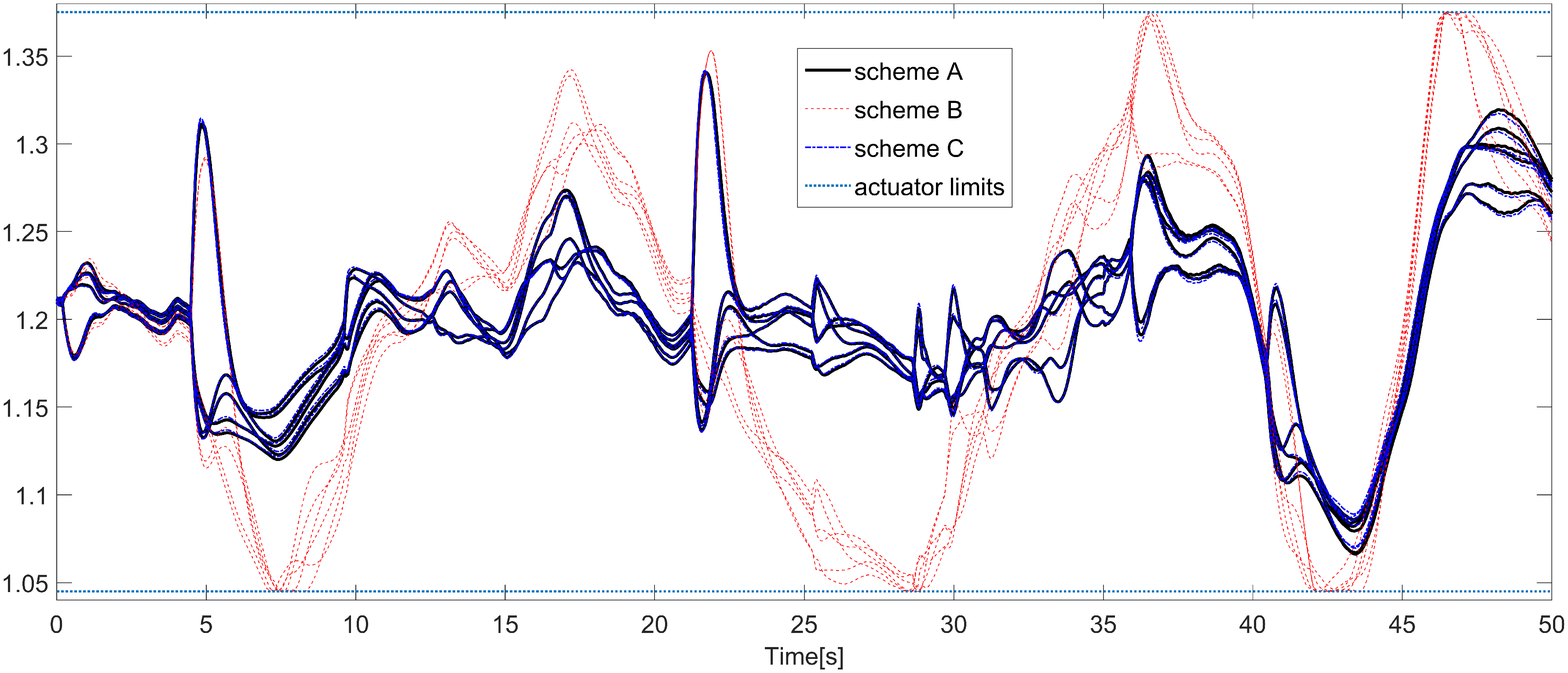}\label{actuators}}
    \caption{(a): Vertical acceleration trajectories using the three MPC schemes for DiM; (b): Hexapod actuators length using the three MPC schemes for DiM. In each simulation, the length of six actuators is drawn in the same color.}
    \label{fs1}
\end{figure*}

\begin{figure*}[!ht]
    \centering
    \subfloat[]{\includegraphics[width=.48\textwidth]{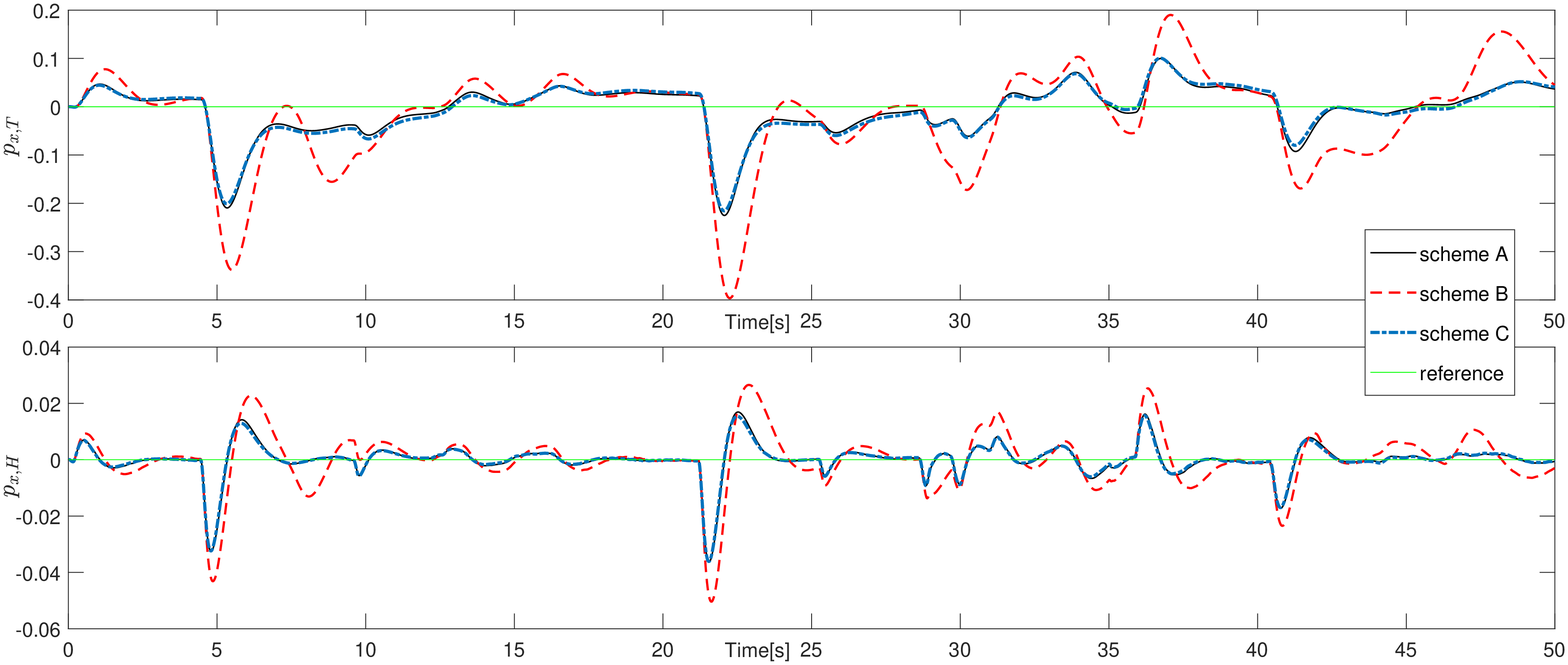}\label{px}}
    \quad
    \subfloat[]{\includegraphics[width=.48\textwidth]{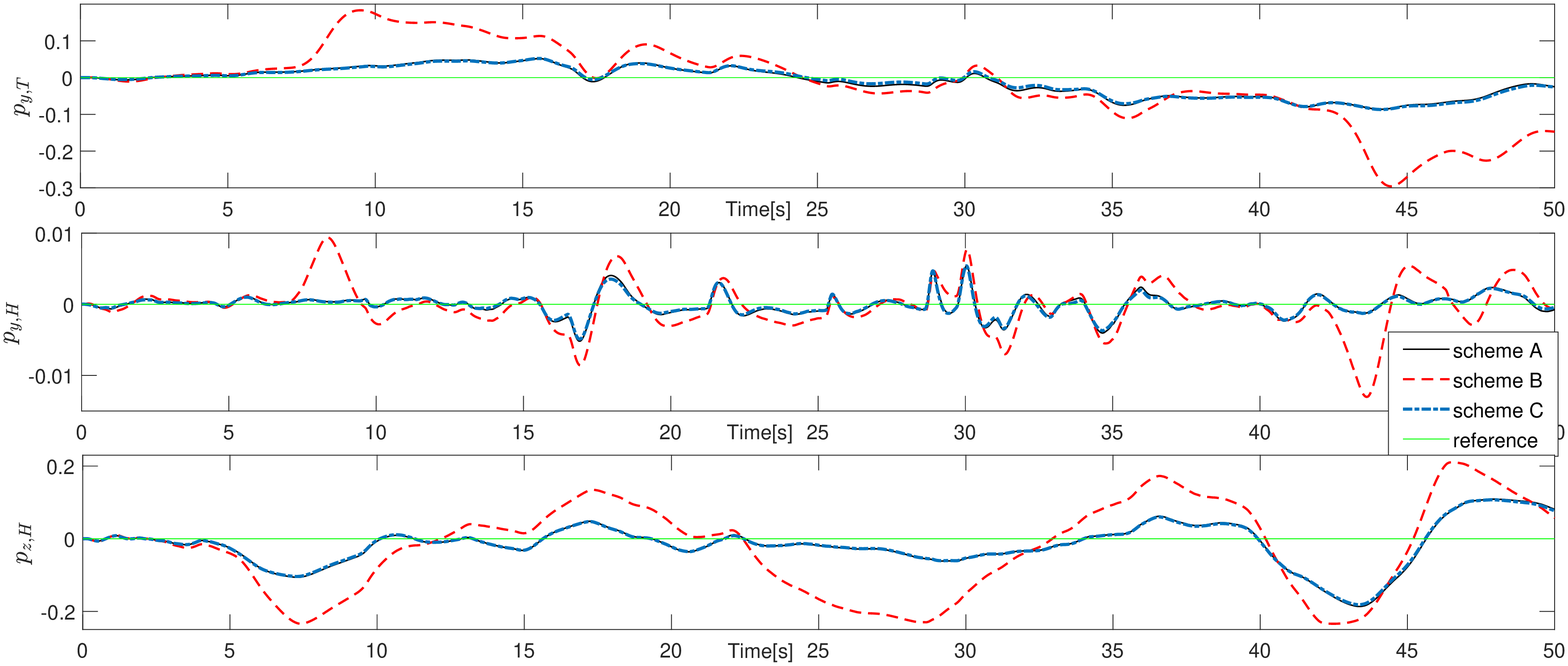}\label{pypz}}
    \caption{(a): Tripod and hexapod longitudinal displacements using the three MPC schemes for DiM; (b): Tripod lateral and hexapod lateral and vertical displacements using the three MPC schemes for DiM}
    \label{fs2}
\end{figure*}

Table \ref{CPT} gives the maximum computation time for one RTI step during the simulation using the three schemes. For the proposed input MB scheme, the computation time for condensing is only slightly higher than that of the nonuniform grid scheme, but still much smaller than that of the standard NMPC. The computation time for solving the QP problem is also much smaller than that of the standard NMPC due to greatly reduced DoFs. Since the input MB scheme considers more state and path constraints, its time for solving QP is slightly higher than that of nonuniform grid scheme.

\begin{table}[!ht]
\centering
\caption{Maximum computation times [ms] per RTI step using the three NMPC schemes for inverted pendulum}
\label{CPT}
\resizebox{.8\linewidth}{!}{%
\begin{tabular}{cccc}
\hline
Scheme  &A &  B & C \\ \hline
Shooting   & 0.26         & 0.26             & 0.26   \\ 
Condensing & 1.8         & 0.05             & 0.37   \\ 
QP solving & 2.26         & 0.09             & 0.15   \\ 
Total      & 4.32         & 0.40             & 0.78   \\ \hline
\end{tabular}
}
\end{table}

It is interesting to compare performance of a specific parameterization of scheme B that shows similar computational time to that of scheme C. This can be achieved by increasing the number of intervals for non-uniform grid NMPC (scheme B) to $M=42$ with
\begin{equation*}
    \begin{aligned}
    I = &[0,1,2,3,4,5,6,7,8,9,10,11,12,13,14,15,16,17,18,\\
    &19,20,21,22,23,24,26,28,32,35,37,40,42,44,46,48,\\
    &50,52,55,60,65,70,75,80]
    \end{aligned}
\end{equation*}
The state and control trajectories are shown in Fig. \ref{traj_compare_new} and the KKT values are shown in Fig. \ref{kkt_compare_new}. It can be observed that, the state and control trajectories of scheme B with $M=42$ become very close to the standard NMPC (scheme A) but the KKT values are still much larger. Therefore, we can conclude that the developed input MB NMPC algorithm has better closed-loop and numerical performance than the non-uniform grid NMPC when the number of input blocks is the same as the number of discretization intervals ($M=10$ for both schemes), whereas if the closed-loop and computational performance are similar, the developed algorithm is able to achieve an improved solution optimality using even less number of input blocks ($M=10$ v.s. $M=42$).

\subsection{DiM}
\begin{figure}[htb]
    \centering
    \includegraphics[width=.5\textwidth]{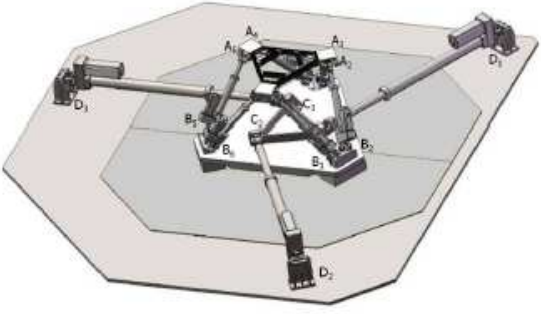}
    \caption{An illustration of the mechanical structure of the DiM platform.}
    \label{fig:DiM}
\end{figure}

The second application is the Motion Cueing algorithm for a nine DoF dynamic driving simulator, which is described in \cite{9DoFDiM}. Dynamic platforms are designed to reproduce in a safe and reproducible environment the driving experience, bridging the virtual and real world through the Human-in-the-Loop (HiL). Their effectiveness is strongly related to the capability of reproducing a realistic sensations of driving with dedicated motion strategies, namely \emph{Motion Cueing Algorithms} (MCAs), that are also deputed to the management of the working area of the device, by keeping its movements within the given operating limits (\emph{Washout Action}). A promising approach to MCA is based on MPC \cite{bruschetta2017fast,9DoFDiM}: to account for the human perception, a model of the vestibular sensory system is considered, together with the mechanical structure of the platform. Perceived accelerations on the vehicle are used as reference for the controller, which is deputed to balance between tracking those and reducing the overall platform motion.
The mechanical structure of DiM is shown in Fig. \ref{fig:DiM}: a hexapod is mounted on a tripod that moves on a flat and stiff surface, sliding on airpads and providing longitudinal, lateral and yaw movements, while the hexapod provides high frequency longitudinal, lateral and yaw movements as well as pitch and roll rotations.
The complete dynamic model consists of a human vestibular model that is characterized by a linear state-space realization and a reference transition model.

The compact dynamic equation of DiM can be summarized as following:
\begin{equation}
\begin{aligned}\label{DiM}
&\dot{\xi}_{VEST}= A \xi_{VEST} + B u_{VEST}\\
&\ddot{p}_D=\ddot{p}_T+R_T \ddot{p}^T_H+2 \omega_T \times (R_T \dot{p}^T_H)  
\end{aligned}
\end{equation}
where the first equation describes the dynamics of the vestibular model, which converts the actual longitudinal, lateral, vertical accelerations and roll, pitch and yaw angular velocities into the corresponding perceived signals in brain. The second equation describes the reference transition model and $p_D$ represents the position of the driver head with respect to a frame fixed to the ground, $p_T$ the position of the centre of the tripod referred to the ground frame, $p_H^T$ the position of the hexapod center with respect to the frame fixed to the tripod, $R_T$ the rotation matrix which describes the frame fixed to the tripod with respect to the ground and  $\omega_T$ its related angular velocity. 

Observe that the longitudinal, lateral acceleration and yaw angular velocity can be written as a sum of low frequency (tripod) component and high frequency (hexapod) one as
\begin{equation}\label{6dof}
    \begin{aligned}
    &a_x=a_{x,T}+a_{x,H},\\
    &a_y= a_{y,T}+a_{y,H},\\
    &\dot{\phi}=\dot{\phi}_T+\dot{\phi}_H.
    \end{aligned}
\end{equation}
We introduce a simplified six DoF model where these two components are obtained by means of a low-pass and high-pass filter. This model has an input $u(\cdot)$ of dimension $6$, a state $x(\cdot)$ of dimension $30+6n_s$, where $n_s$ is the order of the filters, an output $y(\cdot)$, composed of the human perceived signals and the DiM motion states, of dimension $30$. Finally, the MPC problem is formulated as
\begin{equation}
    \begin{aligned}\label{DiM_MPC}
        \min_{x(\cdot),u(\cdot)}  &\int_{t_0}^{t_f} (\norm{y(t)-y_{REF}(t)}_{Q}^2\text{d}t + \norm{y(t_f)-y_{REF}(t_f)}_{Q_f}^2 \\
        s.t. \quad &x(t_0) = \hat{x}_0, \\
            & \text{dynamics given in (\ref{DiM}}),\eqref{6dof}, \\
        & q_-\leq q(x(t))\leq q_+,\\
	  & x_- \leq x(t) \leq x_+,\\
        & u_-\leq u(t)\leq u_+
    \end{aligned}
\end{equation}
where $y_{REF}(t), t \in [t_0,t_f]$, is the reference signal that the DiM must track, $q(x(\cdot))$ is the function which maps the state into the actuators lengths, $q_-,q_+$ the lower and upper actuator length limits, while $x_-,u_-,x_+,u_+$ are the upper and lower bounds for inputs and states specified in  \cite{9DoFDiM}.
For the three schemes we use the following configurations:
\begin{enumerate}
    \item \textit{standard unblocked RTI (scheme A)}: the sampling and shooting interval time is $T_s=10$ ms and $N=50$.
    \item \textit{Non-uniform grid (scheme B)}: the index of nonuniform grid is given by $I = [0, 1, 10, 50]$. Hence the number of stages is $M=3$. The weights for the cost function are based on those used by scheme A, scaled by the length of each interval.
    \item \textit{input MB (scheme C)}: the index of input block is given by the same I as in scheme B. Hence, the number of stages is $N=50$ and the DoFs of NMPC is $M=3$. The weights for the cost function are identical to those used by scheme A.
\end{enumerate}

Fig. \ref{traj_compare_DiM} shows the perceived vertical acceleration using the three NMPC schemes with respect to a given reference trajectory. The performance of scheme A and C are almost identical, while that of scheme B are very bad when actuators reach their limits. In fact, as shown in Fig. \ref{actuators}, scheme A and scheme C are able to maintain the actuators in a more neutral position, while with the scheme B actuators reach their limits in many instants. This can be also verified by Fig. \ref{px} and \ref{pypz}, where scheme A and C can maintain the platform in a more neutral position. On the other hand, scheme B drives the platform in a much larger spaces with worse tracking accuracy, which may harm the mechanical structure of DiM.

\begin{figure}[htb]
    \centering
    \includegraphics[width=\linewidth]{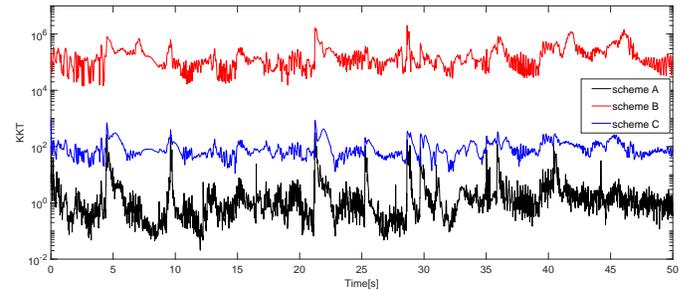}
    \caption{Closed-loop KKT values using the three MPC schemes for problem \eqref{DiM_MPC}} 
    \label{kkt_compare_DiM}
\end{figure}

Fig. \ref{kkt_compare_DiM} shows the KKT values using the three NMPC schemes for DiM. Similar to the inverted pendulum example, the proposed scheme C can maintain the optimality of on-line solution to a large extent by solving a sparse optimization problem with much less DoFs. In Table \ref{CPT-DiM}, the maximum computation time per RTI step using the three NMPC schemes for DiM is given. Both scheme B and C can considerably reduce the computation time, achieve real-time feasibility (sampling time is $10$ ms). However, scheme C has a much superior control performance as demonstrated in previous figures with only a slight higher computation time than scheme B.  

\begin{table}[thb]
\centering
\caption{Maximum computation times [ms] per RTI step using the three MPC schemes for DiM}
\label{CPT-DiM}
\resizebox{.8\linewidth}{!}{%
\begin{tabular}{cccc}
\hline
Scheme  &A &  B & C \\ \hline
Shooting   & 2.91         & 2.91             & 2.91   \\ 
Condensing & 9.65         & 0.08             & 1.22   \\ 
QP solving & 5.81         & 0.17             & 0.27   \\ 
Total      & 18.07        & 3.16             & 4.40   \\ \hline
\end{tabular}
}
\end{table}

\section{Conclusion}\label{sec5}
In this paper, an computationally efficient MB strategy for multiple shooting based NMPC has been proposed. Input MB is introduced during the shooting step, resulting in a multi-stage QP subproblem with reduced DoFs. A tailored condensing algorithm has also been proposed to exploit such sparsity structure, thus reducing the computation complexity of condensing from quadratic to linear in the number of discretization nodes. The convergence and stability property of the proposed scheme is also addressed. Through a detailed comparison with the nonuniform grid scheme, the proposed MB strategy has been shown to achieve better discretization accuracy, constraint fulfillment, easier tuning, and sparser problem structure. 

By using two numerical examples, it has also been shown that the proposed MB strategy is able to significantly reduce on-line computation time while maintaining solution accuracy and control performance. We assume a block structure making the input MB NMPC convergent, recursively feasible and stable. Future developments can focus on an proper choice of input block structure with feasibility and stability guarantees. The error bound on solution and optimality can also be addressed when input MB is applied.

\bibliographystyle{IEEEtran}
\bibliography{myIET}

\begin{thebibliography}{10}

\bibitem{quirynen2015multiple}
Quirynen, R., Vukov, M., Diehl, M.
\newblock `Multiple shooting in a microsecond'.
\newblock In: Multiple Shooting and Time Domain Decomposition Methods.
  (Springer,  2015. pp.~ 183--201

\bibitem{lazutkin2018approach}
Lazutkin, E., Geletu, A., Li, P.: `An approach to determining the number of
  time intervals for solving dynamic optimization problems', \emph{Industrial
  \& Engineering Chemistry Research},  2018, \textbf{57}, (12), pp.~4340--4350

\bibitem{paiva2015adaptive}
Paiva, L.T., Fontes, F.: `Adaptive time-mesh refinement in optimal control
  problems with state constraints', \emph{Discrete and Continuous Dynamical
  Systems},  2015, \textbf{35}, (9), pp.~4553--4572

\bibitem{potena2018non}
Potena, C., Della.Corte, B., Nardi, D., Grisetti, G., Pretto, A.
\newblock `Non-linear model predictive control with adaptive time-mesh
  refinement'.
\newblock In: Simulation, Modeling, and Programming for Autonomous Robots
  (SIMPAR), 2018 IEEE International Conference on. (IEEE,  2018. pp.~ 74--80

\bibitem{lee2018mesh}
Lee, K., Moase, W., Manzie, C.: `Mesh adaptation in direct collocated nonlinear
  model predictive control', \emph{International Journal of Robust and
  Nonlinear Control},  2018, \textbf{28}, (15), pp.~4624--4634

\bibitem{cagienard2007move}
Cagienard, R., Grieder, P., Kerrigan, E.C., Morari, M.: `Move blocking
  strategies in receding horizon control', \emph{Journal of Process Control},
  2007, \textbf{17}, (6), pp.~563--570

\bibitem{gondhalekar2009controlled}
Gondhalekar, R., Imura, J.i., Kashima, K.: `Controlled invariant
  feasibility—a general approach to enforcing strong feasibility in mpc
  applied to move-blocking', \emph{Automatica},  2009, \textbf{45}, (12),
  pp.~2869--2875

\bibitem{gondhalekar2010least}
Gondhalekar, R., Imura, J.i.: `Least-restrictive move-blocking model predictive
  control', \emph{Automatica},  2010, \textbf{46}, (7), pp.~1234--1240

\bibitem{shekhar2015optimal}
Shekhar, R.C., Manzie, C.: `Optimal move blocking strategies for model
  predictive control', \emph{Automatica},  2015, \textbf{61}, pp.~27--34

\bibitem{yu2016stable}
Yu, M., Biegler, L.T.: `A stable and robust nmpc strategy with reduced models
  and nonuniform grids', \emph{IFAC-PapersOnLine},  2016, \textbf{49}, (7),
  pp.~31--36

\bibitem{paiva2018sufficient}
Paiva, L.T., Fontes, F.A.: `A sufficient condition for stability of
  sampled--data model predictive control using adaptive time--mesh refinement',
  \emph{IFAC-PapersOnLine},  2018, \textbf{51}, (20), pp.~104--109

\bibitem{bock1984multiple}
Bock, H.G., Plitt, K.J.: `A multiple shooting algorithm for direct solution of
  optimal control problems', \emph{IFAC Proceedings Volumes},  1984,
  \textbf{17}, (2), pp.~1603--1608

\bibitem{diehl2002real}
Diehl, M., Bock, H.G., Schl{\"o}der, J.P., Findeisen, R., Nagy, Z.,
  Allg{\"o}wer, F.: `Real-time optimization and nonlinear model predictive
  control of processes governed by differential-algebraic equations',
  \emph{Journal of Process Control},  2002, \textbf{12}, (4), pp.~577--585

\bibitem{binder2001introduction}
Binder, T., Blank, L., Bock, H.G., Bulirsch, R., Dahmen, W., Diehl, M., et~al.
\newblock `Introduction to model based optimization of chemical processes on
  moving horizons'.
\newblock In: Online optimization of large scale systems. (Springer,  2001.
  pp.~ 295--339

\bibitem{zanelli2017forces}
Zanelli, A., Domahidi, A., Jerez, J., Morari, M.: `Forces nlp: an efficient
  implementation of interior-point methods for multistage nonlinear nonconvex
  programs', \emph{International Journal of Control},  2017, pp.~ 1--17

\bibitem{stellato2018osqp}
Stellato, B., Banjac, G., Goulart, P., Bemporad, A., Boyd, S.
\newblock `Osqp: An operator splitting solver for quadratic programs'.
\newblock In: 2018 UKACC 12th International Conference on Control (CONTROL). (,
   2018. pp.~ 339--339

\bibitem{frison2013fast}
Frison, G., Jorgensen, J.B.
\newblock `A fast condensing method for solution of linear-quadratic control
  problems'.
\newblock In: Decision and Control (CDC), 2013 IEEE 52nd Annual Conference on.
  (IEEE,  2013. pp.~ 7715--7720

\bibitem{andersson2013general}
Andersson, J.
\newblock `A general-purpose software framework for dynamic optimization'.
\newblock PhD thesis, Arenberg Doctoral School, KU Leuven, Department of
  Electrical Engineering (ESAT/SCD) and Optimization in Engineering Center,
  2013

\bibitem{diehl2001real}
Diehl, M.
\newblock `Real-time optimization for large scale nonlinear processes'.
\newblock Heidelberg University,  2001

\bibitem{diehl2005nominal}
Diehl, M., Findeisen, R., Allg{\"o}wer, F., Bock, H.G., Schl{\"o}der, J.P.:
  `Nominal stability of real-time iteration scheme for nonlinear model
  predictive control', \emph{IEE Proceedings-Control Theory and Applications},
  2005, \textbf{152}, (3), pp.~296--308

\bibitem{mayne2000constrained}
Mayne, D.Q., Rawlings, J.B., Rao, C.V., Scokaert, P.O.: `Constrained model
  predictive control: Stability and optimality', \emph{Automatica},  2000,
  \textbf{36}, (6), pp.~789--814

\bibitem{pannocchia2011conditions}
Pannocchia, G., Rawlings, J.B., Wright, S.J.: `Conditions under which
  suboptimal nonlinear mpc is inherently robust', \emph{Systems \& Control
  Letters},  2011, \textbf{60}, (9), pp.~747--755

\bibitem{graichen2010stability}
Graichen, K., Kugi, A.: `Stability and incremental improvement of suboptimal
  mpc without terminal constraints', \emph{IEEE Transactions on Automatic
  Control},  2010, \textbf{55}, (11), pp.~2576--2580

\bibitem{allan2017inherent}
Allan, D.A., Bates, C.N., Risbeck, M.J., Rawlings, J.B.: `On the inherent
  robustness of optimal and suboptimal nonlinear mpc', \emph{Systems \& Control
  Letters},  2017, \textbf{106}, pp.~68--78

\bibitem{chen2018matmpc}
Chen, Y., Bruschetta, M., Picotti, E., Beghi, A.
\newblock `Matmpc-a matlab based toolbox for real-time nonlinear model
  predictive control'.
\newblock In: European Control Conference. (,  2019.

\bibitem{chen2017fast}
Chen, Y., Cuccato, D., Bruschetta, M., Beghi, A.
\newblock `A fast nonlinear model predictive control strategy for real-time
  motion control of mechanical systems'.
\newblock In: Advanced Intelligent Mechatronics (AIM), 2017 IEEE International
  Conference on. (IEEE,  2017. pp.~ 1780--1785

\bibitem{chen2017inexact}
Chen, Y., Cuccato, D., Bruschetta, M., Beghi, A.
\newblock `An inexact sensitivity updating scheme for fast nonlinear model
  predictive control based on a curvature-like measure of nonlinearity'.
\newblock In: Decision and Control (CDC), 2017 IEEE 56th Annual Conference on.
  (IEEE,  2017. pp.~ 4382--4387

\bibitem{chen2018adaptive}
Chen, Y., Bruschetta, M., Cuccato, D., Beghi, A.: `An adaptive partial
  sensitivity updating scheme for fast nonlinear model predictive control',
  \emph{IEEE Transactions on Automatic Control},  2018,

\bibitem{chen2018efficient}
Chen, Y., Frison, G., van Duijkeren, N., Bruschetta, M., Beghi, A., Diehl, M.:
  `Efficient partial condensing algorithms for nonlinear model predictive
  control with partial sensitivity update', \emph{IFAC-PapersOnLine},  2018,
  \textbf{51}, (20), pp.~406--411

\bibitem{ferreau2014qpoases}
Ferreau, H.J., Kirches, C., Potschka, A., Bock, H.G., Diehl, M.: `qpoases: A
  parametric active-set algorithm for quadratic programming',
  \emph{Mathematical Programming Computation},  2014, \textbf{6}, (4),
  pp.~327--363

\bibitem{quirynen2015autogenerating}
Quirynen, R., Vukov, M., Zanon, M., Diehl, M.: `Autogenerating microsecond
  solvers for nonlinear mpc: A tutorial using acado integrators', \emph{Optimal
  Control Applications and Methods},  2015, \textbf{36}, (5), pp.~685--704

\bibitem{9DoFDiM}
Bruschetta, M., Maran, F., Beghi, A.: `Nonlinear, mpc-based motion cueing
  algorithm for a high-performance, nine-dof dynamic simulator platform',
  \emph{IEEE Transactions on Control Systems Technology},  2016, \textbf{25},
  pp.~1--9

\bibitem{bruschetta2017fast}
Bruschetta, M., Maran, F., Beghi, A.: `A fast implementation of mpc-based
  motion cueing algorithms for mid-size road vehicle motion simulators',
  \emph{Vehicle system dynamics},  2017, \textbf{55}, (6), pp.~802--826

\end{thebibliography}

\end{document}